\documentclass[12pt]{article}

\usepackage{amsmath}
\usepackage{amssymb}
\usepackage{amsfonts}
\usepackage{latexsym}

\newcommand{\e}{{\varepsilon}}

\title{Lower bounds for quasianalytic functions, {\rm II}.
\\ The Bernstein quasianalytic functions.}
\date{}

\author{A. Borichev, F. Nazarov and  M. Sodin\thanks
{Supported by the Israel Science Foundation of the Israel Academy
of Sciences and Humanities under Grant No. 37/00-1.}}

\begin{document}
\maketitle

\begin{abstract}
Let $\mathcal F$ be a class of functions with the uniqueness
property: if $f\in \mathcal F$ vanishes on a set $E$ of positive
measure, then $f$ is the zero function. In many instances, we
would like to have a quantitative version of this property, e.g.
a lower bound for $|f|$ outside a small exceptional
set. Such estimates are well-known and useful for polynomials,
complex- and real-analytic functions, exponential polynomials. In
this work we prove similar results for the Denjoy-Carleman and the
Bernstein classes of quasianalytic functions.

In the first part, we considered quasianalytically smooth functions.
Here, we deal with classes of functions characterized by exponentially
fast approximation by polynomials whose degrees belong to a given very
lacunar sequence. We also prove the polynomial spreading lemma 
and a comparison lemma which are of a certain interest on their own.
\end{abstract}

\section{Introduction and the results}

Let $f$ be a continuous function on $[-1,1]$, and let
$$
E_n(f) = \min_{P\in \mathcal P_n} \|f-P\|_{[-1,1]}
$$
be the {\em approximating sequence of the function} $f$. Here
$\mathcal P_n$ is a space of all algebraic polynomials of degree
$\le n$, and the norm $\|\, . \,\|_F$ denotes the uniform norm
$\|\,.\,\|_{C(F)}$ on the set $F$. A classical result of S.~Bernstein
\cite{Bern1, Bern2} states that if for some $\beta>0$
$$
E_n(f) \le e^{-\beta n} \eqno (1.1)
$$
when $n$ runs through a subsequence $\{n_j\}\subset \mathbb N$, and if the
function $f$ vanishes on a subset of $[-1,1]$ of positive measure,
then $f$ is the zero function. If a sequence $\{n_j\}$ is not too
lacunary:
$$
\limsup_{j\to \infty} \frac{n_{j+1}}{n_j} \le \Delta <\infty,
$$
then condition (1.1) describes a class of real-analytic functions
on $[-1,1]$ with analytic extension into a certain complex
neighbourhood of $[-1,1]$ whose size depends on the quotient
$\beta/\Delta$.

The functions satisfying (1.1) are called the {\em Bernstein
quasianalytic functions}. Having the uniqueness property,
generally speaking, they do not posses any smoothness. They were
studied by Bernstein \cite{Bern1, Bern2}, Beurling
\cite{Beurling}, Mergelyan \cite[Chapter~VII]{Mergelyan},
Ple\'sniak \cite{Ples}, by no means is this list complete.

Here, we give an asymptotic upper bound for the size of the level
sets
$$
m_f(t)=\left| \left\{ x\in [-1,1]: \, |f(x)|\le t \right\}\right|
$$
for $t=E^*_{n_j}(f)$ where $E_n^*(f)=\max\left(E_n(f), e^{-n}\right)$.
The main result follows:

\medskip\par\noindent{\bf Theorem A. }{\em Suppose $f$ is a Bernstein 
quasianalytic function
satisfying condition {\rm (1.1)} with sufficiently lacunar sequence
$\{n_j\}$:
$$
\lim_{j\to \infty} \frac{n_{j+1}}{n_j} = +\infty. \eqno (1.2)
$$
Then
$$
\lim_{j\to \infty} \frac{|\log m_f (E^*_{n_{j+1}}(f))|}
{|\log m_f(E^{*}_{n_j}(f))|} = +\infty. \eqno (1.3)
$$
}

\medskip\par\noindent{\bf Corollary. }{\em In the assumptions of 
Theorem~A, we have
$$
\lim_{j\to \infty} \frac{\log|\log m_f(E^*_{n_j}(f))|}{j} = +\infty.
\eqno(1.4)
$$
}

Mention that relations (1.3) and (1.4) are
fulfilled with $E_{n_j}(f)$ instead of $E^*_{n_j}(f)$ if we assume
additionally that
$$
\lim_{j\to \infty} \frac{|\log E_{n_{j+1}}(f)|}{|\log E_{n_j}(f)|}=+\infty.
$$
This follows from the proof of the Theorem~A given below.

The more lacunary is the sequence $\{n_j\}$ in Theorem~A, the worse is our 
bound (1.4). This is natural since as our second result shows, the 
Bernstein quasianalytic functions may have deep zeros of prescribed 
flatness:

\medskip\par\noindent{\bf Theorem B.}
{\em Given decreasing functions $\varphi,\psi:[1,+\infty)\to(0,+\infty)$,
$$
\lim_{t\to\infty}\varphi(t)=\lim_{t\to\infty}\psi(t)=0,
$$
there exist a function $f\in C[-1,1]$ and a subsequence
$\{n_j\}\subset \mathbb N$ such that for $n\in \{n_j\}$
$$
E_n(f) \le \psi (n), \eqno (1.5)
$$
and
$$
|f(x)| \le e^{-n}, \qquad |x|\le \varphi (n). \eqno (1.6)
$$
}

\medskip In particular, if $\psi(s)=e^{-s}$, we get
$$
m_f(E^*_{n_j}(f)) \ge 2\varphi (n_j), \qquad j=1,2,\,...\,.
$$

\medskip\par\noindent{\sc Acknowledgement. } The authors thank
Eugenia  Malinnikova  and Alex\-ander Volberg for useful
discussions.

\section{The polynomial spreading lemma}

\medskip The key ingredient in the proof of Theorem~A is the
following

\medskip\par\noindent{\bf Polynomial spreading lemma. }{\em
Let $P\in \mathcal P_n$, $\|P\|_{[-1,1]}\le 1$, and let
$$
\delta\ge \delta_0>0, \qquad 0<c_0\le c<1,
\qquad 0<\e <\frac{1-c}{2-c},
$$
be some parameters. Suppose $E\subset [-1,1]$ is a measurable
subset of sufficiently small measure
$$
|E| < \kappa (\delta_0, \e, c_0) \eqno (2.1)
$$
such that
$$
\|P\|_E \le e^{-\delta n},
$$
and suppose that $I$ is an interval, $E\subset I \subset [-1,1]$.
Then
$$
\|P\|_I \le e^{-c\delta n}
$$
provided
$$
|I| \le |E|^{\frac{1}{2-c}+\e}. \eqno (2.2)
$$
}

\medskip
If the set $E$ is itself an interval, then (2.2) can be
significantly improved to $ |I| \le |E|^{c+\e} $ (cf. the end of
the proof of comparison lemma below). This can be regarded as a
polynomial version of the Hadamard three circle theorem. An
``ideal statement'' would be
$$
\|P\|_I \le \|P\|^c_E \|P\|^{1-c}_{[-1,1]}
$$
provided that $ E \subset I \subset [-1,1] $ are intervals such
that
$$
|I| \le |E|^c 2^{1-c}\,.
$$
This is too good to be true. Our result gives a reasonable
approximation to such a logarithmic convexity.

The exponent $ \frac1{2-c}+\epsilon $ in (2.2) is larger than the exponent 
$ c+\epsilon $ we need. However, they are close to each other when $ c $ 
is close to one. To obtain Theorem~A we first use a 
dyadic decomposition with a simple stopping-time rule,  
and then apply iteratively the spreading lemma for $c=1-\rho$
with small $\rho>0$ to get the comparison lemma (Section~3) which claims that
under natural conditions
\begin{gather*}
|\{x\in [-1,1]: |P(x)|\le e^{-t\delta n} \|P\|_{[-1,1]}\}|\\
\qquad \qquad \ge |\{x\in [-1,1]: |P(x)|\le e^{-\delta n} 
\|P\|_{[-1,1]}\}|^{t+\gamma}.
\end{gather*}
The proof of Theorem~A is then completed in Section~4.

\medskip\par\noindent{\em Proof of the spreading lemma:}
We use an argument adopted from Nadirashvili's work
\cite{Nadirashvili, Nadirashvili1}. Let $ \eta=|I|/|E| $, and let
$x_0$ be the centre of the interval $I$. Fix $k\ge 0$, and 
consider the Taylor polynomial
$$
P_k(x) = \sum_{j=0}^k \frac{P^{(j)}(x_0)}{j!} (x-x_0)^j, 
$$
and the remainder $ R_k = P - P_k $. Applying to $P_k$ the
classical Remez inequality \cite{Remez} (cf. \cite{BG, DR}), we
get
\begin{eqnarray*}
\|P\|_I &\le& \|P_k\|_I + \|R_k\|_I \\ \\
&\le & (4\eta)^k \|P_k\|_E + \|R_k\|_{[-1,1]} \\ \\
&\le& (4\eta)^k \|P\|_E + \left( 1+ (4\eta)^k\right)
\|R_k\|_{[-1,1]} \\ \\
&\le& (4\eta)^k \left( e^{-\delta n} + 2
\|R_k\|_{[-1,1]}\right).
\end{eqnarray*}
Using the Lagrange formula for the remainder, we have
$$
\|R_k\|_{[-1,1]} \le \left(\frac{|I|}{2} \right)^{k+1}
\frac{\|P^{(k+1)}\|_{[-1,1]}}{(k+1)!} < \left( \frac{e}{2}\,
\frac{|I|}{k+1} \right)^{k+1} \|P^{(k+1)}\|_{[-1,1]}.
$$
Recalling the classical V.~Markov inequality \cite{W.M} for the $(k+1)$-st
derivative of the polynomial $P$ of degree $n$ 
\footnote{There are several relatively simple proofs of this inequality, 
see e.g. \cite{DS} for one of them. In fact, we could use a slightly 
cruder version of Markov's estimare given in \cite[Chapter~VI, Lemma 
4.III]{M} with a proof found by Th.~Bang. }
$$
\|P^{(k+1)}\|_{[-1,1]} \le \frac1{2} \left( \frac{2}{k+1}
\right)^{k+1} n^{2k+2} \|P\|_{[-1,1]},
$$
we get
$$
\|R_k\|_{[-1,1]} \le \frac1{2} \left( \frac{e|I|n^2}{(k+1)^2}
\right)^{k+1},
$$
and then
$$
\|P\|_I \le (4\eta)^k \left\{ e^{-\delta n} +
\left(\frac{e|I|n^2}{(k+1)^2}\right)^{k+1} \right\}. \eqno (2.3)
$$

Now, our requirements to the choice of $k$ are the following:
$$
(4\eta)^k \le e^{\delta (1-c) n} \eqno (2.4)
$$
and
$$
\left( \frac{e|I|n^2}{(k+1)^2}\right)^{k+1} \le e^{-\delta n}.
\eqno (2.5)
$$
Naturally, relations (2.3)--(2.5) yield that
$$
\|P\|_I \le 2 e^{-c\delta n}. \eqno (2.6)
$$

Applying (2.6) to $P^M$ with $M\in \mathbb N$, we obtain
$$
\|P^M\|_I \le 2 e^{-c\delta nM}
$$
or
$$
\|P\|_I=\|P^M\|^{1/M}_I\le 2^{1/M}e^{-c\delta n}=(1+o(1))e^{-c\delta n},
\qquad M\to \infty,
$$
completing the proof of the lemma.

It remains to verify that there exists $k$ satisfying (2.4)--(2.5).
Suppose that for some positive $\lambda$,
$$
\lambda\le \delta \frac{1-c}{\log 4\eta},\eqno (2.7)
$$
$$
\lambda\log\frac{\lambda^2}{e|I|}\ge \delta, \eqno (2.8)
$$
$$
\lambda^2 \ge \frac{|I|}e, \eqno (2.9)
$$
and choose $k$ such that
$$
k\le \lambda n< k+1.
$$
Then (2.4) follows from (2.7):
$$
(4\eta)^k \le (4\eta)^{\lambda n} \le e^{\delta (1-c)n},
$$
and (2.5) follows from (2.8) because the left-hand side of (2.8)
increases as a function of $\lambda$ satisfying (2.9):
\begin{eqnarray*}
\left( \frac{e|I|n^2}{(k+1)^2} \right)^{k+1} &=& \exp\left[ -(k+1)
\log \frac{(k+1)^2}{e|I|n^2} \right] \\ \\
&\le& \exp\left[
-\lambda n \log \frac{\lambda^2}{e|I|} \right] \le \exp[-\delta
n]\,.
\end{eqnarray*}

Without loss of generality, assume that $\delta\le 1$, otherwise we just 
increase the degree $n$ in the statement of the lemma. 
We denote $A=\delta^2/(e|I|)$, $ B = A/\log^2 A $, and set $
\lambda = \delta/ \log B $. We have to check that inequalities
(2.7)--(2.9) hold for this choice of $\lambda$ under the condition that the
length of $E$ (and therefore that of $I$) is sufficiently small,
that is the value $A$ is sufficiently large.

Estimate (2.7) says that
$$
(4\eta)^{1/(1-c)}\le e^{\delta/\lambda} = B =
\frac{\delta^2}{e|I|\log^2\frac{\delta^2}{e|I|}},
$$
or, equivalently,
$$
\left( 4 \frac{|I|}{|E|}\right)^{1/(1-c)} \le \frac{\delta^2}{e|I|
\log^2 \frac{\delta^2}{e|I|} },
$$
$$
|I|^{\frac{2-c}{1-c}}\log^2\frac{\delta^2}{e|I|} \le 
\frac{\delta^2}{4^{1/(1-c)}e}|E|^{1/(1-c)}
$$
that follows from (2.1) and (2.2).

Inequality (2.8) becomes
$$
\frac\delta{\log B } \log \frac {A}{\log^2 B} \ge\delta,
$$
that is
$$
\frac A{\log^2 B} \ge \frac A{\log^2 A},
$$
which is evidently true.

At last, inequality (2.9) becomes
$$
\frac{1}{\log^2 B} \ge \frac{1}{e^2 A}
$$
which is true when $A$ is sufficiently large. \hfill $ \Box
$

\section{Comparison lemma}

Here, we give a corollary to the spreading
lemma which will be needed for the proof of Theorem~A. For
$P\in\mathcal P_n$, we set
$$
\mathcal E_P(\delta) = \left\{x\in [-1,1]: |P(x)|\le e^{-\delta n
} \|P\|_{[-1,1]} \right\}.
$$

\medskip\par\noindent{\bf Comparison Lemma. }{\em
Let $P\in \mathcal P_n$ and let
$$
\delta \ge \delta_0>0, \qquad 0 < t_0 \le t
<1, \qquad  0<\gamma < 1-t
$$
be some parameters. Suppose that the length of the set $\mathcal E_P(\delta)$
is
\newline\noindent
$(\delta_0, \gamma, t_0)$-sufficiently small. Then
$$
|\mathcal E_P(t\delta)| \ge |\mathcal E_P(\delta)|^{t+\gamma}. \eqno (3.1)
$$
}

\smallskip\par\noindent{\em Proof: } Without loss of generality we assume 
that
$\|P\|_{[-1,1]}=1$. First of all, we prove a weaker result:

\medskip\par\noindent{\bf Claim. }{\em If
$$
0<c_0\le c<1,\qquad 0<\e<\frac {1-c}{2-c},
$$
and if
$$
|\mathcal E_P(\delta)| \le \kappa(\delta_0, \e, c_0),
$$
where $\kappa$ has the same value as in the spreading lemma, then
$$
|\mathcal E_P(c\delta)| \ge |\mathcal E_P(\delta)|^{\frac 1{2-c}+\e}.
$$
}

Set $E=\mathcal E_P(\delta)$ and choose
an integer
$$
N \ge \kappa^{-1}(\delta_0, \e, c_0).
$$
Let
$\mathcal J$ be the collection of all maximal $N$-adic subintervals
$I$ of $[-1,1]$ such that
$$
|E\cap I|^{\frac 1{2-c}+\e} \ge |I|.
$$
Then the ``remainder set''
$$
F=E\setminus \bigcup_{I\in \mathcal J} (E\cap I)
$$
has zero length. Indeed, for any $\xi>0$ we can cover $F$ by
disjoint $N$-adic intervals $J_\alpha$ of length $|J_\alpha| < \xi$, and
$$
|F\cap J_\alpha| \le |J_\alpha|\cdot \xi^{(\frac1{2-c}+\e)^{-1}-1}.
$$
Summing up by $\alpha$, we see that
$$
|F|\le \xi^{(\frac1{2-c}+\e)^{-1}-1},\qquad \xi>0,
$$
and hence, $|F|=0$.

Now, we apply the spreading lemma to the sets $E\cap I$, $I\in
\mathcal J$. Since $I_0=[-1,1] \notin \mathcal J$, we have
$$
|E\cap I|\le |I| \le \frac1{N} \le \kappa (\delta_0, \e_0, c_0),
$$
and the conditions of the lemma are satisfied. Hence,
$$
\|P\|_I \le e^{-c\delta n}, \qquad I\in \mathcal J,
$$
that is
$$
\bigcup_{I\in \mathcal J} I \subset \mathcal E_p(c\delta).
$$

Since $I$ is the maximal interval, its $N$-adic ``supinterval''
$I^*$ does not belong to $\mathcal J$, that is
$$
|E\cap I|^{\frac 1{2-c}+\e} \le |E\cap I^*|^{\frac 1{2-c}+\e}
\le|I^*|=N|I|,
$$
or $|I|\ge N^{-1}|E\cap I|^{\frac 1{2-c}+\e}$. Therefore,
\begin{eqnarray*}
|\mathcal E_P(c\delta) | &\ge & \sum_{I\in \mathcal J} |I| \ge
\frac1{N} \sum_{I\in \mathcal J} |E\cap I|^{\frac 1{2-c}+\e} \\ \\
&\ge & \frac1{N} \left( \sum_{I\in \mathcal J} |E\cap I|
\right)^{\frac 1{2-c}+\e} = \frac{|E|^{\frac 1{2-c}+\e}}{N} =
\frac{|\mathcal E_P(\delta)|^{\frac 1{2-c}+\e}}{N}\,.
\end{eqnarray*}
Increasing slightly $\e$, we get the claim.

\medskip
Now we choose an integer $M$ and $\e>0$ in such a way that
$$
\Bigl(\frac 1{2-t^{1/M}}+\e\Bigr)^M\le t+\gamma,
$$
set $c_0=t_0$, $c=t^{1/M}$, and apply the claim $M$ times. We get
$$
|\mathcal E_P(t\delta)| \ge |\mathcal E_P(\delta)|^
{\bigl(\frac 1{2-t^{1/M}}+\e\bigr)^M}
\ge |\mathcal E_P(\delta)|^{t+\gamma},
$$
unless
$$
|\mathcal E_P(t\delta)| \ge \kappa(\delta_0, \e, c_0).
$$
In both cases we get (3.1), and thus the lemma is proved. \hfill $ \Box $

\section{Proof of Theorem~A}

Theorem~A is a simple corollary to the comparison lemma.

Let $P_{n}$ be the polynomials of the best approximation to $f$, that is
$$
\|f-P_n\|_{[-1,1]}=E_n(f).
$$
Without loss of generality, we assume that $\|f\|_{[-1,1]}=1$, and
$$
1/2 \le \|P_{n_j}\|_{[-1,1]} \le 2.
$$ 
Now, according to the definition of
$E^*_n(f)$ and (1.1),
$$
E^*_{n_j}(f)=e^{-\delta_j n_j}, \qquad \min(\beta,1) \le \delta_j\le 1.
$$
Due to the lacunarity condition (1.2), for any $\e>0$,
$$
\frac{1}{4} E^*_{n_{j-1}}(f) \ge {(4 E^*_{n_j}(f))}^\e, \qquad
j\ge j(\e).
$$
Then
\begin{gather*}
\left\{x\in[-1,1]\!:|f(x)|\le E^*_{n_j}(f) \right\}\\ 
\qquad \qquad \qquad \qquad \qquad \subset
\left\{x\in[-1,1]\!:|P_{n_j}(x)|\le 4E^*_{n_j}(f) \|P_{n_j}\|_{[-1,1]} 
\right\}
\end{gather*}
and
\begin{gather*}
\left\{x\in[-1,1]\!:|f(x)|\le E^*_{n_{j-1}}(f) \right\}\\ 
\qquad \qquad \qquad \qquad \supset
\left\{x\in[-1,1]\!:|P_{n_j}(x)|\le 
\frac{1}{4} E^*_{n_{j-1}}(f) \|P_{n_j}\|_{[-1,1]}
\right\}.
\end{gather*}
Therefore, applying the comparison lemma to the
polynomials $P_{n_j}$ with $t=\gamma=\e$, we get for sufficiently
large $j$:
\begin{gather*}
m_f(E^*_{n_{j-1}}(f)) \ge 
\left|\left\{x\in[-1,1]\!:|P_{n_j}(x)|\le 
\frac{1}{4} E^*_{n_{j-1}}(f) \|P_{n_j}\|_{[-1,1]}\right\} \right| \\ \\
\qquad \qquad \ge \left|\left\{x\in[-1,1]\!:|P_{n_j}(x)|\le 4E^*_{n_j}(f) 
\|P_{n_j}\|_{[-1,1]}
\right\}\right|^{2\e} \ge m_f^{2\e}( E^*_{n_j}(f)).
\end{gather*}
This proves the theorem. \hfill $\Box $

\section{Proof of Theorem~B}

We start with

\medskip\par\noindent{\bf Lemma.}
{\em Let $Q$ be a polynomial, $Q(0)=0$. Then for any odd
positive integer $n$ and any sufficiently large integer $l\ge
l_0(n)$ there is a polynomial $P$ of degree at most $ln\deg Q$
such that
$$
\|P\|_{[-1,1]} \le C_1\cdot\frac{\|Q\|_*}{n},
$$
and
$$
|(Q+P)(t)| \le C_1\cdot\frac{(2n |t|)^{l+1}}{(l+1)!} \|Q\|_*\,, \qquad
|t|\le \frac 1n.
$$
}

\smallskip Here $\|Q\|_*$ means the sum of the absolute values of the
coefficients of $Q$, and $C_1$ is a constant.

\smallskip\par\noindent{\em Proof: } First, we prove a special case of the 
lemma with
$Q(t)=t$. Set
$$
\Phi_n(w)=n\sin\left( \frac{1}{n}\arcsin w\right)\,,
\qquad  w=u+iv\,.
$$
The functions $\Phi_n$ are analytic in the unit disc, continuous up to its
boundary and uniformly bounded. Furthermore,
$$
|\Phi_n(w)| \lesssim |w|\,, \qquad |w|\le 1\,,
$$
where the notation $A\lesssim B$ means that
$A\le C \cdot B$ for a positive numerical constant $C$.
Also,
$$
|\Phi_n^{(k)}(w)| \lesssim  2^k\,,
\qquad |w| \le \frac{1}{2}\,, \quad k\in {\mathbb Z}_+.
$$
Set
$$
\Phi_{n,l}(w) = \sum_{k=0}^l \frac{\Phi_n^{(k)}(0)}{k!} w^k\,.
$$
By Abel's theorem, the polynomials $\Phi_{n,l}(u)$ converge to $\Phi_n(u)$
uniformly in $u\in [-1,1]$, so that
$$
|\Phi_{n,l}(u)| \lesssim 1\,,
\qquad u\in [-1,1]\,, \quad l\ge l_0(n)\,,
\eqno (5.1)
$$
and
$$
|(\Phi_n - \Phi_{n,l})(u)| \le \frac{|u|^{l+1}}{(l+1)!}
\|\Phi_n^{l+1}\|_{[-u,u]} \lesssim \frac{(2|u|)^{l+1}}{(l+1)!},
\qquad |u|\le \frac{1}{2}\,. \eqno (5.2)
$$

Let $t=\sin\left( \frac{1}{n} \arcsin u\right)$, then
$u=u_n(t)=\sin(n\arcsin t)$. Since $n$ is odd, $u_n(t)$ is a polynomial of
degree $n$, and
$$
|u_n(t)|\le\min(1,n|t|),\qquad |t|\le 1\,. \eqno (5.3)
$$
Indeed, it is sufficient to verify that
$$
n\sin t-\sin(nt)\ge 0,\qquad 0\le \sin t\le \frac 1n.
$$
This inequality holds for $t=0$, and taking the derivative
of the left-hand side we get
$$
n[\cos t-\cos (nt)]
$$
which is non-negative on the interval $[0,\arcsin\frac 1n]$.

Set
$$
R_{n,l}(t) = -\frac{1}{n} \left( \Phi_{n,l}\circ u_n\right)(t).
$$
This is a polynomial of degree at most $ln$. We have
\begin{eqnarray*}
\left|t+R_{n,l}(t)\right| &=& \left| \frac1{n} \left( \Phi_n -
\Phi_{n,l}\right)(u_n(t)) \right| \\ \\ &\stackrel{(5.2)}\lesssim&
\frac{(2|u_n(t)|)^{l+1}}{(l+1)!} \stackrel{(5.3)}\lesssim
\frac{(2n|t|)^{l+1}}{(l+1)!}\,, \qquad \qquad  |t| \le \frac{1}{n},
\end{eqnarray*}
and
$$
|R_{n,l}(t)| \stackrel{(5.1)}\lesssim \frac1{n}, \qquad |t|\le 1\,.
$$
This proves the special case of the lemma.

The general case follows if we set
$$
Q(t)=\sum_{j\ge 1} c_jt^j,
$$
and
$$
P(t)= \sum_{j\ge 1} c_j R_{n,l}(t^j).
$$
\hfill $ \Box $

\medskip\par\noindent{\bf Corollary. }{\em Given 
a polynomial $Q$, $Q(0)=0$, $\e>0$, $N<\infty$,
and given a function $\varphi$ decreasing to
zero at infinity, there exist $M>N$,
and a polynomial $P$, $P(0)=0$,
such that $\deg P\le M$,
$$
\|P\|_{[-1,1]} \le \e,
$$
and
$$
|(Q+P)(t)| \le e^{-2M}\,,
\qquad |t|\le \varphi(M).
$$
}

\smallskip\par\noindent{\em Proof:} We apply the lemma with
$$
n=\frac{C_1\cdot\|Q\|_*}{\e},
$$
and $l$ such that
\begin{gather*}
M=ln\deg Q>N,\\
\varphi(M)\le\frac 1{2n},\\
(l+1)!\ge  C_1\cdot e^{2M}\|Q\|_*.
\end{gather*}
We get a polynomial $P$ of degree at most $M$ such that
$$
\|P\|_{[-1,1]} \le \e,
$$
and
$$
|(Q+P)(t)|\le\frac{C_1\|Q\|_*}{(l+1)!}\le e^{-2M}\,, 
\qquad |t| \le \varphi(M)\le\frac 1{2n}.
$$
This establishes the corollary. \hfill $ \Box $

\smallskip\par\noindent{\em Proof of Theorem~B: }
We use the corollary in an inductive procedure. We build
the function $f$ in the form
$$
f=\sum_{j\ge 1} P_j,
$$
where $P_j$ are polynomials such that for an increasing sequence
of integers $\{n_j\}$ we have
$\deg P_j\le n_j$,
$\|P_j\|_{[-1,1]} < \psi(n_{j-1})/2$,
and
$$
\Bigl|\Bigl(\sum_{1\le j \le m} P_j\Bigr)(x)\Bigr|\le
e^{-2n_m}\,, \qquad |x|\le \varphi(n_m).
$$

We start with $P_1(x)=x$. After $P_1$,\ldots, $P_m$ have been
chosen, set
\begin{gather*}
Q=\sum_{1\le j \le m} P_j,\\
\e=\psi(n_m)/2,\\
N=n_{m},
\end{gather*}
and  get $n_{m+1}=M>n_m$ and $P_{m+1}$ from the corollary.

We may always assume that
$$
\psi(x)\le e^{-2x}\,,\qquad x\ge 1\,,\eqno (5.5)
$$
otherwise, from the very beginning, we replace $\psi (x)$ by $\min
(e^{-2x},\psi(x))$. Therefore,
$$
\sum_{j\ge m} \psi(n_{j-1})/2\le \psi(n_m)\,.\eqno (5.6)
$$

Now, we check that $f$ and $\{n_j\}$ satisfy conditions (1.5) and (1.6) of
Theorem~B. For $n=n_m$ we have
$$
\sum_{j\ge m+1} \|P_j\|_{[-1,1]} \le \sum_{j\ge m+1} \psi(n_{j-1})/2
\stackrel{(5.6)}\le \psi(n)
$$
that proves (1.5). Finally, for $|x|\le \varphi (n)$
$$
|f(x)| \le\Bigl|\Bigl(\sum_{1\le j \le m} P_j\Bigr)(x)\Bigr|+
\sum_{j\ge m+1} \|P_j\|_{[-1,1]}\le e^{-2n}+\psi(n)\le e^{-n}
$$
that proves (1.6). \hfill $ \Box $

\section{Remarks and questions}

{\bf 6.1 Beurling's theorem.} Beurling \cite[p.~396--403]{Beurling} gave
a general quasianalyticity condition which contains those of
Bernstein and Denjoy-Carleman. Here, we formulate a special case
of his result. Given a sequence $1\ge e_n\downarrow 0$, consider
the Bernstein class
$$
\mathcal F_{\{e_n\}} = \{f\in C[-1,1]:\, E_n(f)\le e_n \}
\setminus \{0\}
$$

\medskip\par\noindent{\bf Theorem }(Beurling). {\em The class
$\mathcal F_{\{e_n\}}$ contains no function vanishing on a
subset of positive measure if and only if
$$
\sum_{n\ge 1} \frac{\log^-e_n}{n^2} = +\infty
$$
where $\log^-a = \max(\log\frac1a,0)$.}

Beurling's proof uses the Laplace transform combined with the harmonic
estimation in the "if part" and the Paley-Wiener theorem in the "only
if part". One can extract a quantitative estimate from his
proof which however is essentially weaker than Theorem~A above and
Theorem~B from the part I \cite{NSV} of this work.

It seems to be interesting to obtain another proof of Beurling's theorem
by means of the constructive function theory and to get its quantitative
version in a sharp form.

\medskip\par\noindent{\bf 6.2 Potential theory approach. } A minute
reflection suggests that there could be a natural generalization of the
spreading lemma and comparison lemma for the logarithmic potentials of 
probability measures.

Let $ u $ be a subharmonic function in the complex plane $ \mathbb C $ 
with compactly supported Riesz measure $ \mu $, $ \mu(\mathbb C) \le 1 $, 
and let $ u\Big\vert_{[0,1]} \le 0 $. Let $ E\subset [0,1] $ be a subset 
of positive measure such that $ u\Big\vert_E \le -\delta $.

\smallskip\par\noindent{\bf Problem. }{\em
Given $ c $, $ 0<c<1 $, estimate from below the length of the set 
$ \{x\in [0,1]: \, u(x)\le -c\delta \} $. 
}

\smallskip
An ``ideal'' lower bound would be $ |E|^c $, which is true in 
the trivial limiting cases $ c=0 $ and $ c=1 $. Probably, our 
polynomial comparison lemma (combined with a suitable atomization of the 
Riesz measure $ \mu $) yields an ``asymptotic'' lower bound $ 
|E|^{c+\epsilon} $. However, it seems more natural to treat this problem 
by means of potential theory.  

One can also ask a similar question replacing the unit interval by the 
unit disk. In this case, probably, one should deal with capacity instead 
of linear measure.

\bigskip

\medskip\par\noindent (A.B.:) Laboratoire de Math\'ematiques Pures,

\par\noindent Universit\'e Bordeaux I,

\par\noindent Talence, 33405,

\par\noindent France

\par\noindent {\textit{\small E-mail}: \texttt{\small
borichev@math.u-bordeaux.fr}}

\medskip\par\noindent (F.N.:) Department of Mathematics,

\par\noindent Michigan State
University,

\par\noindent East Lansing, MI 48824,

\par\noindent U.S.A.

\par\noindent {\textit{\small E-mail}: \texttt{\small
fedja@math.msu.edu}}

\medskip
\par\noindent (M.S.:) School of Mathematical Sciences,

\par\noindent Tel Aviv University,

\par\noindent Ramat Aviv, 69978,

\par\noindent Israel

\par\noindent {\textit{\small E-mail}: 
\texttt{\small sodin@post.tau.ac.il}}

\bigskip


\begin{thebibliography}{2}
\bigskip

\bibitem{Bern1} {\sc S. Bernstein},
\newblock{\em Sur la d\'efinition et les propri\'et\'es des fonctions
analytiques d'une variable r\'eelle},
\newblock{Math. Ann. {\bf 75} (1914), 449--468.}

\bibitem{Bern2} {\sc S. Bernstein},
\newblock{\em Le\c{c}ons sur les propri\'et\'e extr\'emales et la
meilleure approximation des fonctions analytiques d'une variable
r\'eelle}, \newblock{Paris, Gauthier Villars, 1926.}

\bibitem{Beurling} {\sc A. Beurling},
\newblock{\em Quasianalyticity. Mittag-Leffler Lectures on Complex
Analysis}, (1977-1978)
\newblock{Collected Works, Vol \rm {I}, Birkh\"auser, Boston, 1989.}

\bibitem{BG} {\sc Yu. Brudnyi and M. Ganzburg},
\newblock{\em An extremal problem for polynomials of $n$ variables},
\newblock{Izv. AN SSSR (Mat) {\bf 37} (1973), 344--355. (Russian)}

\bibitem{DR} {\sc R. M. Dudley and B. Randol},
\newblock{\em Implications of pointwise bounds on polynomials},
\newblock{Duke Math. J. {\bf 29} (1962), 455--458.}

\bibitem{DS} {\sc R. J. Duffin and A. C.  Schaeffer}, 
\newblock{\em A refinement of an inequality of the brothers Markoff}, 
\newblock{Trans. Amer. Math. Soc. {\bf 50} (1941), 517--528.}

\bibitem{M} {\sc Sz. Mandelbrojt},
\newblock{\em S\'eéries adh\'erentes, r\'egularisation des suites,
applications },
\newblock{Gauthier-Villars, Paris, 1952.}

\bibitem{W.M} {\sc V. Markov},
\newblock{\em On functions least deviated from zero on a given interval},
\newblock{St.-Petersburg, 1892 (Russian)}.
\newblock{German transl. Math. Ann. {\bf 77} (1916), 213--258.}

\bibitem{Mergelyan} {\sc S. Mergelyan},
\newblock{\em Some questions of the constructive function theory},
\newblock{Proc. of the Steklov Math. Inst., Vol. \rm{XXXVII}, Moscow,
1951. (in Russian)}

%\bibitem{MMR} {\sc G. V. Milovanovi\'c, D.~S.~Mitrinovi\'c and 
%Th.~M.~Rassias},
%\newblock{\em Topics in polynomials: extremal problems, inequalities, 
%zeros}, 
%\newblock{World Scientific, Singapore, 1994.}

\bibitem{Nadirashvili} {\sc N. Nadirashvili},
\newblock{\em On a generalization of Hadamard's three-circle theorem},
\newblock{Vestnik Mosk. Un-ta. (Mat.) {\bf 31} no. 3 (1976),
39--42.}

\bibitem{Nadirashvili1} {\sc N. Nadirashvili},
\newblock{\em Estimating solutions of elliptic equations with analytic
coefficients bounded on some sets},
\newblock{Vestnik Mosk. Un-ta. (Mat.) {\bf 34} no. 2 (1979),
42--46.}

\bibitem{NSV} {\sc F. Nazarov, M. Sodin and A. Volberg},
\newblock{\em Lower bounds for quasianalytic functions, I. How to
control smooth functions?}
\newblock{arXiv CA/0208238.}

\bibitem{Ples} {\sc W. Ple\'sniak},
\newblock{\em Quasianalytic functions in the sense of Bernstein},
\newblock{Dissert. Math. {\bf 147} (1977). }

\bibitem{Remez} {\sc E. J. Remez},
\newblock{\em Sur une propri\'ete des polyn\^omes de Tschebycheff},
\newblock{Commun. Inst. Sci. Kharkov {\bf 13} (1936), 93--95.}


\end{thebibliography}
\end{document}